\documentclass{amsart}

\usepackage{graphicx, hyperref}

\newtheorem{theorem}{Theorem}[section]
\newtheorem{lemma}[theorem]{Lemma}
\newtheorem{proposition}[theorem]{Proposition}

\newcommand{\thistheoremname}{}
\newtheorem*{genericthm}
{\thistheoremname}
\newenvironment{namedtheorem}[1]
  {\renewcommand{\thistheoremname}{#1}%
   \begin{genericthm}}
  {\end{genericthm}}

\theoremstyle{definition}

\theoremstyle{remark}
\newtheorem{remark}[theorem]{Remark}

\numberwithin{equation}{section}

\newcommand{\tr}{{\hbox{\rm tr}}}
\newcommand{\sn}{{\hbox{\rm sn}}}
\newcommand{\disc}{{\hbox{\rm disc}}}
\newcommand{\im}{{\hbox{\rm im}}}
\newcommand{\Norm}{{\hbox{\rm Norm}}}
\newcommand{\Hil}{{\hbox{\rm Hil}}}
\newcommand{\U}{{\hbox{\rm U}}}
\newcommand{\grO}{{\hbox{\rm O}}}
\newcommand{\SO}{{\hbox{\rm SO}}}

\begin{document}

\title{On the Spinor norm in the unitary groups}

\author{Ngo Van Dinh}

\address{Department of Mathematics and Informatics, Thai Nguyen University of Sciences, Thai Nguyen, Vietnam}

\email{dinh.ngo@tnus.edu.vn}





\begin{abstract}
Let $F$ be a field of odd characteristic, $E$ be a finite extension of $F$ equipped an involution with subfield of fixed points $E_0$ containing $F$ and $V$ be a finite dimensional $E$-vector space with a non-degenerate hermitian form $h$. We show a link between the spinor norm in the unitary group $\U(V,h)$ and the calculus of determinants and discriminants. Then we show a formula which links the spinor norm in $\U(V,h)$ and the spinor norm in the orthogonal group $\grO(V,b_h)$ defined by a non-degenerate symmetric bilinear form $b_h$ associated to $h$.
\end{abstract}

\maketitle

\section{Introduction}
Let $F$ be a field of characteristic not equal to $2$ and $E$ be a finite extension of $F$, equipped an involution $\bar{ }$ which fixes all the elements of $F$. Denote by $E_0$ the subfield of fixed points of $E$ by the involution. We fix a non-zero $F$-linear form $\mu_0$ from $E_0$ to $F$ and put $\mu=\mu_0\circ \tr_{E/E_0}$, where $\tr_{E/E_0}$ is the trace form from $E$ to $E_0$.

Let $V$ be a finite dimesional vector space over $E$ and $h: V\times V\to E$ be a non-degenerate hermitian form on $V$. Considering $V$ as an $F$-vector space, we have an associated non-degenerate symmetric bilinear form $b_h$ defined by
\[b_h(x,y)=\mu(h(x,y)), \text{ for all } (x,y)\in V\times V.\]
Denote $\U(V,h)$ the unitary group defined by the hermitian form $h$ and $\grO(V,b_h)$ the orthogonal group defined by the symmetric form $b_h$. Then every element $\sigma \in \U(V,h)$ is clearly also an element of $\grO(V,b_h)$.

Let $\sn$ be the spinor norm in the group $\grO(V,b_h)$. It is a homomorphism from $\grO(V,b_h)$ to $F^\times/(F^\times)^2$. The spinor norm has important applications \cite{Kneser}. Zassenhaus \cite{Zassenhaus} found a direct definition of $\sn$ which links this norm with the calculus of determinants and discriminants: for each $\tilde{\sigma}\in \grO(V,b_h)$, we have
\begin{equation}\label{Zassenhausformula}
\sn(\tilde{\sigma})=\begin{cases}\disc(b_h) \text{ if } \tilde{\sigma}=-1, \\
\disc\left(b_h|_{\underset{i\geq 1}{\bigcup}\ker (1+\tilde{\sigma})^i}\right) {\det}_F \left(\frac{1+\tilde{\sigma}}{2}|_{\underset{i\geq 1}{\bigcap}\im (1+\tilde{\sigma})^i}\right) \text{ otherwise},\end{cases}
\end{equation}
where $\disc(b_h)$ and $\disc\left(b_h|_{\underset{i\geq 1}{\bigcup}\ker (1+\tilde{\sigma})^i}\right)$ are respectively the discriminant of $b_h$ and the discriminant of the restriction of $b_h$ over subspace $\underset{i\geq 1}{\bigcup}\ker (1+\tilde{\sigma})^i$.

It's well known that there is an anti-hernitian form $h'$ on $V$ such that the unitary groups $\U(V,h')$ and $\U(V,h)$ coincide \cite{Dieudonne}. We summarize here Wall's construction \cite{Wall} of a spinor norm on the unitary group $\U(V,h')$. Let $\sigma$ be a non-trivial element of $\U(V,h)$. Denote $V_{\sigma}$ the image of the transformation $1-\sigma$. If $V_\sigma$ is an $E$-vector subspace of dimension $r$ then we say that $\sigma$ is an \emph{element of dimension $r$}. By definition, each element of dimension $1$, denoted $s_{(v,\varphi)}$, is defined by
\[s_{(v,\varphi)}(x)=x-\varphi h'(v,x)v, \text{ for all } x\in V,\]
where $v$ is a non-zero element of the space $V_{s_{(v,\varphi)}}$ and $\varphi$ is an element of $E^{\times}$ such that $\varphi^{-1}-\bar{\varphi}^{-1}=h'(v,v)$. For each element $\sigma\in \U(V,h')$ of dimension $r>0$, let $f_\sigma$ be the sesquilinear form (with respect to the involution $\bar{ }$ ) on $V_\sigma$ defined by
\[\begin{split}
f_\sigma: V_\sigma \times V_\sigma &\to E,\\
\left(x-\sigma(x),y-\sigma(y)\right) &\mapsto h'(x-\sigma(x),y).
\end{split}\]
Then $\sigma$ can be written as a product of one-dimensional elements \cite[Lemma 3]{Wall},
\[\sigma=s_{(v_1,\varphi_1)}s_{(v_2,\varphi_2)}\ldots s_{(v_r,\varphi_r)},\]
where the vectors $v_1,v_2,..., v_r$ form a orthogonal basis of $V_\sigma$ with respect to the form $f_\sigma$ anf $v_1$ can be chosen as any non-isotropic vector of $V_\sigma$. A such decomposition of $\sigma$ is called \emph{Cayley decomposition} of $\sigma$.

Let $a$ be a fixed vector of $V$ and suppose $\sigma \in \U(V,h')$ is an element of dimension $r>0$ with a Cayley decomposition
\[\sigma=s_{(v_1,\varphi_1)}s_{(v_2,\varphi_2)}\ldots s_{(v_r,\varphi_r)},\]
where the vector $v_i, i=1,2,...,r,$ is chosen such that $h'(a,v_i)$ is either $0$ or $1$. Then the class $\varphi_1\varphi_2...\varphi_rE_0^{\times}$ in $E^{\times}/E_0^{\times}$ depends only on $\sigma$ \cite[Lemma 5]{Wall}. This class is called \emph{spinor norm of $\sigma$}, denoted by $\sn_E(\sigma)$:
\begin{equation}
\sn_E(\sigma)=\varphi_1\ldots \varphi_r E_0^\times.
\end{equation}

The spinor norm of the identity of $\U(V,h')$ is defined to be $E_0^{\times}$. Then we have a homomorphism $\sn_E$ of $\U(V,h')$ in $E^{\times}/E_0^{\times}$ \cite[Lemma 6, Lemma 7]{Wall}, called \emph{spinor norm} in $\U(V,h')$.

The goal of this work is to compare the spinor norm in the orthogonal group $\O(V,b_h)$ and the spinor norm in the unitary group $\U(V,h)$. A natural way to do so is to get a formula similar to Zassenhau's formula (\ref{Zassenhausformula}).

\begin{namedtheorem}{Proposition}
Let $s=s_{(v,\varphi)}$ be an one-dimensional element of $\U(V,h')$. Then we have
\begin{equation}\label{FormulaTheorem1}
\disc \left(h'\vert_{\underset{n\geq 1}{\bigcup}\ker(1+s)^n}\right){\det}_E\left(\frac{1+s}{2}\vert_{\underset{n\geq 1}{\bigcap}\im(1+s)^n}\right)=\varphi \mod(E_0^{\times}).
\end{equation}
\end{namedtheorem}

\begin{namedtheorem}{Theorem}
For all $\sigma\in \U(V,h')$, we have 
\begin{equation}\label{FormulaTheorem2}
\sn(\tilde{\sigma})=\Norm_{E/F}(\sn_E(\sigma)),
\end{equation}
where $\Norm_{E/F}$ is the homomorphism of $E^{\times}/E_0^{\times}$ in $F^\times/(F^\times)^2$ induced by the norm of $E$ over $F$.
\end{namedtheorem}

Note that, for the case where $E$ is a quadratic extension of $F$, this observation has been given in \cite[Chapter 10, Theorem 1.5]{Scharlau} with an error proof. We give here another proof for this link in general case.

This paper is based on the research which is part of the doctoral dissertation \cite{Ngo} of the author. The results are useful in the study of supercuspidal representations of spin groups over a $p$-adic field, where some calculations arise involving the restriction of the spinor norm to unitary groups contained in the orthogonal group under study. The author is grateful to Corinne Blondel for her support, advice and interest in this work at various times.

\section{Proof of the proposition}
In order to prove the formula (\ref{FormulaTheorem1}), we distinguish two cases: $v$ is an isotropic vector and $v$ is not one with respect to the form $h'$. In the first case, the element $s=s_{(v,\varphi)}$ is called a \emph{transvection} of $V$. We have $\varphi^{-1}-\bar{\varphi}^{-1}=h'(v,v)=0$, then the spinor norm $\sn_E$ is trivial at $s$.

Now we calculate the left side of (\ref{FormulaTheorem1}). Let $x$ be an element of $\ker(1+s)$. Then $x$ belongs to the one-dimensional $E$-vector subspace of $V$ generated by $v$, \textit{i.e.}, $x=kv$ for some $k\in E$. Since $v$ is isotropic, we have $x+x=0$ and then $x=0$. This follows that the subspace $\underset{n\geq 1}{\bigcup}\ker\left(1+s\right)^n$ is trivial. Furthermore, there exists a basis $\{v,v',w_1,...,w_{n-2}\}$ of $V$ such that $v'$ is also isotropic with respect to $h'$, $h'(v,v')=1$ and $h'(v,w_i)=0, i=1,...,n-2$. Calculate the determinant in this basis, we have
\[{\det}_E\left(\frac{1+s}{2}\right)=\det\begin{pmatrix}
1&-\frac{\varphi}{2}\\ 0&1
\end{pmatrix}.\]
By the calculating above, we see that the formula (\ref{FormulaTheorem1}) holds for the first case.

In the second case, we have an orthogonal decomposition $V=(Ev)\perp (Ev)^{\perp}$ where $(Ev)^{\perp}$ is the orthogonal complement of the line $Ev$ in $V$ with respect to $h'$. Let $x$ be an element of $\ker(1+s)$. Then $x=kv$ for some $k\in E$ and we have $2kv-\varphi h'(v,kv)v=0$. Since $\varphi^{-1}-\bar{\varphi}^{-1}=h'(v,v)$, we have
\[k[1+\varphi\bar{\varphi}^{-1}]v=0.\]
It follows that the subspace $\underset{n\geq 1}{\bigcup}\ker(1+s)^n$ is either zero or the line $Ev$. In the first situation, we have 
\[{\det}_E\left(\frac{1+s}{2}\right)=\frac{1+\varphi\bar{\varphi}^{-1}}{2}\]
while, in the second, we have $\varphi\bar{\varphi}^{-1}=-1$ and
\[\disc\left(h'|_{Ev}\right)=\left(\varphi^{-1}-\bar{\varphi}^{-1}\right)\mod (E_0^{\times}).\]
Then it is easy to see that the formula (\ref{FormulaTheorem1}) holds in these two situations.

\section{Proof of the theorem}
Firstly we note that in the case where $\varphi\bar{\varphi}^{-1}=-1$ the element $s=s_{(v,\varphi)}$ is the reflection of $V$ defined by vector $v$ with respect to $h'$, \textit{i.e.}, it is the linear transformation of $V$ such that $s(v)=-v$ and $s(x)=x$ for all $x\in V$ such that $h'(v,x)=0$. If this is not the case then we have the subspace $\underset{n\geq 1}{\bigcup}\ker(1+s)^n$ is zero as the discussion in the proof of the Proposition. Then by (\ref{Zassenhausformula}) we have 
\[\sn(s)={\det}_F\left(\frac{1+s}{2}\right)\mod (F^{\times})^2\]
and by (\ref{FormulaTheorem1}) we have
\[\sn_E(s)={\det}_E\left(\frac{1+s}{2}\right)\mod (E_0^{\times}).\]
This give us the following lemma which is the affirmation of the formula (\ref{FormulaTheorem2}) for the one-dimensional element which are not reflections.

\begin{lemma}\label{one-dimensionalcase}
Let $s=s_{(v,\varphi)}$ be an one-dimensional element of $\U(V,h')$ which is not a reflection. Then
\[\sn(s)=\Norm_{E/F}\left(\sn_E(s)\right).\]
\end{lemma}

\subsection{Quadratic case} We consider now the case where $E_0=F$. In this case, $\sn$ and $\Norm_{E/F}\circ\sn_E$ are two homomorphisms of $\U(V,h')$ in $F^{\times}/(F^{\times})^2$. Then it suffices to verify the equality (\ref{FormulaTheorem2}) for the one-dimensional elements which generate $\U(V,h')$.

By the Lemma \ref{one-dimensionalcase}, we only need to verify the formula for the reflections of $V$. Let $u=s_{(v,\varphi)}$ be a reflection of $V$, \textit{i.e.}, $\varphi\bar{\varphi}^{-1}=-1$. Then the identity (\ref{FormulaTheorem2}) become
\[\disc(b_h|_{Ev})=\Norm_{E/F}(\disc(h'|_{Ev}))\mod (F^{\times})^2.\]
Note that, in this case, we have $h'=\delta h$ where $\delta\in E$ and $\{1,\delta\}$ forms a orthogonal $F$-basis of $E$ with respect to $b_h$. Identifying the space $Ev$ with $E$, the restriction of $h$ to this space is a hermitian form on $E$. Then we have $h(x,y)=ax\bar{y}, \forall x,y\in E$, for some $a\in E_0$. This follows that $\disc(h'|_E)=a\delta\mod (E_0^{\times})$ and
\[\disc(b_h|_E)={\det}_F\begin{pmatrix}
2a&0\\ 0&-2a\delta^2
\end{pmatrix}\mod (F^{\times})^2=\Norm_{E/F}(a\delta)\mod (F^{\times})^2.\]
That means the formula (\ref{FormulaTheorem2}) holds for the quadratic case:

\begin{proposition}\label{quadraticcase}
If $E_0=F$ then we have
\[\sn(u)=\Norm_{E/F}(\sn_E(u)), \text{ for all } u\in \U(V,h').\]
\end{proposition}

\begin{remark}
Let $u$ be an element of $\U(V,h')$. Suppose ${\det}_E(u)=\alpha$. Then  $\Norm_{E/E_0}(\alpha)=1$. By the Hilbert's theorem 90, there exists an unique element $\beta\in E^{\times}$ up to a scalar in $E_0$ such that $\alpha=\beta\bar{\beta}^{-1}$. Then we have a homomorphism
\[\Hil: \Norm_{E/E_0}^{-1}(1) \to E^{\times}/E_0^{\times}, \alpha\mapsto \beta E_0^{\times}, \text{ where } \alpha=\beta\bar{\beta}^{-1}.\]
The spinor norm $\sn_E$ is in fact the composition of the determinant and the homomorphism $\Hil$, \textit{i.e.}, we have
\[\sn_E(u)=\Hil({\det}_E(u)), \forall u\in U(V,h').\]
In order to prove this identity we only need to verify it for the one-dimensional elements of $\U(V,h')$. Let $u=s_{(v,\varphi)}$ be an one-dimensional element of $\U(V,h')$. If $v$ is an isotropic vector then the identity is evident since ${\det}_E(u)=1$ and the spinor norm of $u$ is trivial by definition. If $v$ is non-isotropic then we have
\[{\det}_E(u)=1-\varphi h'(v,v)=1-\varphi (\varphi^{-1}-\bar{\varphi}^{-1})=\varphi\bar{\varphi}^{-1}.\]
It follows that $\Hil({\det}_E(u))=\varphi E^{\times}_0=\sn_E(u).$ By this point of view of the spinor norm in $\U(V,h')$, we can see that the Proposition \ref{quadraticcase} is similar to \cite[Chapter 10, Theorem 1.5]{Scharlau}. However the proof in \textit{loc.cit.} is not correct since $\tilde{\sigma}\neq \alpha_\beta$ in its notations.
\end{remark}

\subsection{General case} We prove now the Theorem 2 for the general case. For all $x,y\in V$, put
\[h_0(x,y)=\tr_{E/E_0}(h(x,y)).\]
Then $h_0$ is a non-degenerate symmetric bilinear form on $E_0$-vector space $V$. Denote $\SO(V,h_0)$ the group of the rotations of $V$ with respect to $h_0$ and $\sn_{E_0}$ the spinor norm in $\SO(V,h_0)$. Note that we have
\[\U(V,h)\subset \SO(V,h_0) \subset \grO(V,b_h)\]
and, by the Proposition \ref{quadraticcase}, we have
\[\sn_{E_0}(u)=\Norm_{E/E_0}(\sn_E(u)), \forall u\in \U(V,h).\]

For the passage from $E_0$ to $F$, we use the transfer properties of Witt's ring of quadratic spaces \cite[Chapter 9, \S 5]{Scharlau}: Consider $E_0$ as an $E_0$-vector space and denote $\phi_0$ the symmetric bilinear form on $E_0$ defined by
\[\phi_0(x,y)=xy,\forall x,y\in E_0.\]
Then $\mu_0\circ \phi_0$ is a symmetric bilinear form on $F$-vector space $E_0$. Put $\zeta=\disc(\mu_0\circ \phi_0)$. Let $\phi$ be a symmetric bilinear form on an $E_0$-vector space $W$ of dimension $n$. Then $\mu_0\circ \phi$ is also a symmetric bilinear form on $F$-vector space $W$. In this situation, we have \cite[Chapter 9, Theorem 5.12]{Scharlau}
\[\disc(\mu_0\circ \phi)=(\zeta)^n\Norm_{E_0/F}(\disc(\phi)).\]

Return to our situation, let $s$ be a reflection of $E_0$-vector space $V$ with respect to $h_0$. Using the Zassenhaus' formula (\ref{Zassenhausformula}) and the transfer property above we have
\[\sn(s)=\zeta \Norm_{E_0/F}(\sn_{E_0}(s)).\]
Since the reflections generate the group $\SO(V,h_0)$ \cite{Dieudonne}, we obtain
\[\sn(u)=\Norm_{E_0/F}(\sn_{E_0}(u)), \text{ for all } u \in \SO(V,h_0).\]
This completes the proof of the Theorem.

\bibliographystyle{amsplain}

\end{document}